\input amstex
\magnification=1200
\documentstyle{amsppt}
\NoRunningHeads
\define\Diff{\operatorname{Diff}}
\define\Rot{\operatorname{Rot}}
\define\PSL{\operatorname{PSL}}
\define\HS{\operatorname{\Cal H\Cal S}}
\define\TC{\operatorname{\Cal T\Cal C}}
\define\sltwo{\operatorname{\frak s\frak l}(2,\Bbb C)}
\define\ad{\operatorname{ad}}
\define\Hilb{\operatorname{Hilb}}
\define\form{\operatorname{form}}
\define\Pontr{\operatorname{Pontr}}
\define\vir{\operatorname{vir}}
\define\Vect{\operatorname{Vect}}
\define\sgn{\operatorname{sgn}}
\define\End{\operatorname{End}}
\define\bnd{\Cal B}
\define\FR{\operatorname{\Cal F\Cal R}}
\define\K{\Cal K}
\define\Ker{\operatorname{Ker}}
\define\tr{\operatorname{Tr}}
\define\im{\operatorname{Im}}
\define\Vir{\operatorname{Vir}}

\topmatter
\title
Infinite dimensional geometry of $M_1\!=\!\Diff_+(\Bbb S^1)/\PSL(2,\Bbb
R)$\linebreak and $q_R$--conformal symmetries
\endtitle
\author Denis V. Juriev
\endauthor
\date June, 24, 1988\qquad\qquad\qquad\qquad math.RT/9806140
\enddate
\endtopmatter
\document
This article being addressed to specialists in representation theory,
quantum al\-geb\-ra, infinite dimensional differential (K\"ahler) geometry,
mathematical physics and informatics is devoted to various objects of
infinite dimensional geometry related to $q_R$--conformal symmetries, which
analysis being initiated by the author primarily for goals of informatics of
the interactive videosystems acquired mo\-ti\-vations from pure mathematics
later. On the other hand the constructions con\-si\-dered below are originated
ideologically in a certain extent in such branches of mo\-dern mathematical
physics as conformal and string theories and these origins pe\-net\-rated
deeply into the view on nature of the investigated objects. Precisely, the
interrelations of the infinite dimensional geometry of the homogeneous
K\"ahler ma\-ni\-fold $M=\Diff_+(\Bbb S^1)/\Bbb S^1$ [1,2] (which plays a
crucial role in the geometric for\-mu\-lation of conformal and string field
theories [3-5]), its quotient $M_1=\Diff_+(\Bbb S^1)/\PSL(2,\Bbb R)$ and the
approximate representations of the Virasoro algebra by
the $q_R$--conformal symmetries [6] are discussed. I believe that the obtained
results may strengthen and clarify the connections between infinite
dimensional geometry, representation theory, modern mathematical physics
and informatics of interactive systems once suggested in [2,4:III].

\head 1. Preliminary definitions\endhead

\subhead 1.1. The Lie algebra $\sltwo$ and Verma modules over it\endsubhead
Lie algebra $\sltwo$ is a three-dimensional space of $2\times2$ complex
matrices with zero trace supplied by the standard commutator $[X,Y]=XY-YX$,
where the right hand side multiplication is the standard matrix multiplication.
In the basis
$$ l_{-1}=\left(\matrix0&1\\0&0\endmatrix\right),\quad
l_0=\left(\matrix\tfrac12&0\\0&-\tfrac12\endmatrix\right),\quad
l_1=\left(\matrix0&0\\-1&0\endmatrix\right) $$
the commutation relations are of the form:
$$[l_i,l_j]=(i-j)l_{i+j}\qquad\qquad(i,j=-1,0,1).$$
Lie algebra $\sltwo$ is $\Bbb Z$--graded:
$\deg(l_i)=-\ad(l_0)l_i=i$, where $\ad(X)$ is the adjoint action operator
in the Lie algebra: $\ad(X)Y=[X,Y]$. Therefore, $\Bbb Z$--graded modules
over $\sltwo$ are $l_0$--diagonal. A vector $v$ in a $\Bbb Z$--graded
module over the Lie algebra $\sltwo$ is called extremal iff $l_1v=0$
and the linear span of vectors $l_{-1}^nv$ ($n\in\Bbb Z_+$) coincides
with the module itself (i.e. $v$ is a cyclic vector). A $\Bbb Z$--graded
module with an extremal vector (in this case it is defined up to a
multiplier) is called extremal [7]. An extremal module is called the
Verma module [8] iff the action of $l_{-1}$ is free in it, i.e. the vectors
$l_{-1}^nv$ are linearly independent. In the case of the Lie algebra $\sltwo$
the Verma modules are just the same as infinite dimensional extremal
modules. An extremal weight of the Verma module is the number defined
by the equality $l_0v=hv$, where $v$ is the extremal vector. The Verma
modules are defined for all complex numbers $h$ and are pairwise nonisomorphic.
Below we shall consider the Verma modules with real extremal weights only.

The Verma module $V_h$ over the Lie algebra $\sltwo$ with the extremal weight
$h$ may be realized in the space $\Bbb C[z]$ of polynomials of a complex
variable $z$. The formulas for the generators of the Lie algebra $\sltwo$ are
of the form:
$$l_{-1}=z,\quad l_0=z\partial_z+h,\quad l_1=z\partial_z^2+2h\partial_z,$$
here $\partial_z=\frac{d}{dz}$.

The Verma module is nondegenerate (i.e. does not contain any proper
sub\-mo\-dule) iff $h\ne-\tfrac{n}2$ ($n\in\Bbb Z_+$). The Verma module $V_h$
is called unitarizable (or hermitean) iff it admits a structure of the
pre-Hilbert space such that $l_i^*=l_{-i}$. The completion of the
unitarizable Verma module will be denoted by $V^{\Hilb}_h$. The Lie algebra
$\sltwo$ acts in $V^{\Hilb}_h$ by the unbounded operators. Also it is rather
useful to consider the formal Verma modules $V^{\form}_h$, which are realized
in the space $\Bbb C[[z]]$ of formal power series of a complex variable $z$,
whereas the formulas for generators of the Lie algebra $\sltwo$ coincide
with ones above. Note that $V_h\subseteq V^{\Hilb}_h\subseteq V^{\form}_h$
and modules $V_h$, $V^{\Hilb}_h$, $V^{\form}_h$ form the Gelfand triple or
the Dirac equipment of the Hilbert space $V^{\Hilb}_h$. An action of the real
form of the Lie algebra $\sltwo$ generated by the anti-Hermitean operators
$il_0$, $l_1-l_{-1}$ and $i(l_1+l_{-1})$ in the
Hilbert space $V^{\Hilb}_h$ by the unbounded operators is exponentiated to
a unitary representation of the corresponding simply connected Lie group.

In the nonunitarizable Verma module over the Lie algebra $\sltwo$ there
exists the unique (up to a scalar multiple) indefinite sesquilinear
form $(\cdot,\cdot)$ such that $(l_iv_1,v_2)=(v_1,l_{-i}v_2)$ for any two
vectors $v_1$ and $v_2$ from the Verma module. If this sesquilinear form
is nondegenerate (in this case the Verma module is nondegenerate) then it
has a signature $(n,\infty)$, where $n$ is finite, and therefore, there is
defined a Pontryagin completion [9] of the Verma module. The corresponding
module in which the Lie algebra $\sltwo$ acts by the unbounded operators
will be denoted by $V^{\Pontr}_h$. The following chain of inclusions holds:
$V_h\subseteq V^{\Pontr}_h\subseteq V^{\form}_h$. An action of the real
form of the Lie algebra $\sltwo$ generated by the anti-Hermitean (with
respect to the nondegenerate indefinite sesquilinear form $(\cdot,\cdot)$)
operators $il_0$, $l_1-l_{-1}$ and $i(l_1+l_{-1})$ in the
Pontryagin space $V^{\Pontr}_h$ by the unbounded operators is exponentiated
to a pseudounitary representation of the corresponding simply connected
Lie group.

\subhead 1.2. Hidden symmetries in Verma modules over the Lie algebra
$\sltwo$: Lobachevski{\v\i}-Berezin $C^*$--algebra and $q_R$--conformal
symmetries\endsubhead

\proclaim{Proposition 1 [10]}\it In the nondegenerate Verma module $V_h$
over the Lie algebra $\sltwo$ there are uniquely defined the operators
$D$ and $F$ such that
$$\aligned [D,l_{-1}]=1,\quad [D,l_0]=D,\quad [D,l_1]=D^2,\\
[l_{-1},F]=1,\quad [l_0,F]=F,\quad [l_1,F]=F^2.\endaligned$$
If the Verma modules are realized in the space $\Bbb C[z]$ of polynomials
of a complex variable $z$ then
$$D=\partial_z,\qquad F=z\tfrac1{\xi+2h},$$
where $\xi=z\partial_z$. The operators $F$ and $D$ obey the following
relations:
$$[FD,DF]=0,\qquad [D,F]=q_R(1-DF)(1-FD),$$
where $q_R=\tfrac1{2h-1}$. In the unitarizable Verma module ($q_R\!\ne\!0$)
the operators $F$ and $D$ are bounded and $F^*=D$, $D^*=F$.
\endproclaim

\rm
The algebra generated by the variables $t$ and $t^*$ with the relations
$[tt^*,t^*t]=0$ and $[t,t^*]=q_R(1-tt^*)(1-t^*t)$ being the Berezin
quantization of the Lobachevski{\v\i} plane realized in the unit complex disc
(the Poincar\`e realization) [11] is called the Lobachevski{\v\i}-Berezin
algebra. Proposition 1 allows to consider the Lobachevski{\v\i}-Berezin
algebra as a $C^*$--algebra. The Lobachevski{\v\i}-Berezin $C^*$--algebra
was re\-dis\-co\-vered recently by S.Klimek and A.Lesnievsky [12].

\proclaim{Proposition 2 [10]}\it In the nongenerate Verma module $V_h$
over the Lie algebra $\sltwo$ there are uniquely defined the operators $L_n$
($n\in\Bbb Z$) such that
$$[l_i,L_n]=(i-n)L_{i+n}\qquad (i=-1,0,1;\ n\in\Bbb Z),$$
moreover, $L_i=l_i$ ($i=-1,0,1$). If the Verma modules are realized in the
space $\Bbb C[z]$ of polynomials of a complex variable $z$ then
$$L_k=(xi+(k+1)h)\partial_z^k\ (k\ge0),\quad
L_{-k}=z^k\tfrac{\xi+(k+1)h}{(\xi+2h)\ldots(\xi+2h+k-1)}\ (k\ge1),$$
where $\xi=z\partial_z$. The operators $L_n$ obey the following relations:
$$[L_n,L_m]=(n-m)L_{n+m},\text{\rm\ if}\ n,m\ge-1\text{\rm\ or}\
n,m\le1.$$
In the unitarizable Verma module the operators $L_n$ are unbounded and
$L_i^*=L_{-1}$.
\endproclaim

\rm
The operators $L_n$ are called the $q_R$--conformal symmetries. They may
be {\it sym\-bo\-li\-cal\-ly\/} represented in the form:
$$L_n=D^{nh}L_0D^{n(1-h)},\qquad L_{-n}=F^{n(1-h)}L_0F^{nh}.$$
To supply the symbolical recording by a sense one should use the general
com\-mu\-tation relations
$$[L_n,f(D)]=(-D)^{n+1}f'(D)\ (n\ge-1),\qquad
[L_{-n},f(F)]=F^{n+1}f'(D)\ (n\ge-1)$$
for $n=0$.

The commutation relations for the operators $D$, $F$ and the generators of
$q_R$--conformal symmetries with the generators of the Lie algebra $\sltwo$
mean that the families $J_k$ and $L_k$ ($k\in\Bbb Z$), where $J_i=D^i$,
$J_{-i}=F^i$ ($i\in\Bbb Z_+$), are families of tensor operators [9,13] for
the Lie algebra $\sltwo$.

\subhead 1.3. Infinite dimensional $\Bbb Z$--graded Lie algebras:
the Witt algebra $\frak w^{\Bbb C}$ of all Laurent polynomial vector
fields on a circle and the Virasoro algebra $\vir^{\Bbb C}$, its
one-dimensional nontrivial central extension\endsubhead
The Lie algebra $\Vect(S^1)$ is realized in the space of $C^\infty$--smooth
vector fields $v(t)\partial_t$ on a circle $S^1\simeq\Bbb R/2\pi\Bbb Z$
with the commutator
$$[v_1(t)\partial_t,v_2(t)\partial_t]=
(v_1(t)v_2'(t)-v_1'(t)v_2(t))\partial_t.$$
In the basis
$$s_n=\sin(nt)\partial_t,\qquad c_n=\cos(nt)\partial_t,\qquad h=\partial_t$$
the commutation relations have the form:
$$\aligned
[s_n,s_m]&=\tfrac12((m-n)s_{n+m}+\sgn(n-m)(n+m)s_{|n-m|}),\\
[c_n,c_m]&=\tfrac12((n-m)s_{n+m}+\sgn(n-m)(n+m)s_{|n-m|}),\\
[s_n,c_m]&=\tfrac12((m-n)c_{n+m}-(m+n)c_{|n-m|})-n\delta_{nm}h,\\
[h,s_n]&=nc_n,\quad [h,c_n]=-ns_n.
\endaligned$$
Let us denote by $\Vect^{\Bbb C}(S^1)$ the complexification of the Lie algebra
$\Vect(S^1)$. In the basis $e_n=ie^{ikt}\partial_t$ the commutation relations
in the Lie algebra Li $\Vect^{\Bbb C}(S^1)$ have the form:
$$[e_j,e_k]=(j-k)e_{j+k}.$$

It is rather convenient to consider an imbedding of the circle $S^1$ into
the complex plane $\Bbb C$ with the coordinate $z$, so that $z=e^{it}$ on the
circle and the elements of the basis $e_k$ ($k\in\Bbb Z$) are represented by
the Laurent polynomial vector fields: $$e_k=z^{k+1}\partial_z.$$
The $\Bbb Z$--graded Lie algebra generated by the Laurent polynomial vector
fields (i.e. by the finite linear combinations of elements of the basis
$e_k$) is called the Witt algebra and is denoted by $\frak w^{\Bbb C}$.
The Witt algebra $\frak w^{\Bbb C}$ is the complexification of the
subalgebra $\frak w$ of the algebra $\Vect(S^1)$ generated by the
trigonometric polynomial vector fields on a circle $S^1$, i.e. by the finite
linear combinations of elements of the basis $s_n$, $c_n$ and $h$.

The Lie algebra $\Vect(S^1)$ admits a nontrivial one-dimensional central
extension defined by the Gelfand-Fuchs 2-cocycle [14]:
$$c(v_1(t)\partial_t,v_2(t)\partial_t)=
\int_0^{2\pi}(v_1'(t)v_2''(t)-v_2'(t)v_1''(t))\,dt.$$
This extension being continued to the complexification $\Vect^{\Bbb C}(S^1)$
of the Lie algebra $\Vect(S^1)$ and reduced to the subalgebra
$\frak w^{\Bbb C}$ defines a central one-dimensional ex\-tension of the Witt
algebra, which is called the Virasoro algebra and is denoted by
$\vir^{\Bbb C}$. The Virasoro algebra is generated by the elements
$e_k$ ($k\in\Bbb Z$) and the central element $c$ with the commutation
relations:
$$[e_j,e_k]=(j-k)e_{j+k}+\frac{j^3-j}{12}c\delta(j+k)$$
and is the complexification of a central extension $\vir$ of the Lie algebra
$\frak w$.
Because the 2-cocycle
$$\int_0^{2\pi}(v_1(t)v_2'(t)-v_2(t)v_1'(t))\,dt$$
is trivial, the cocycle, which defines the Virasoro algebra, is indeed
equivalent to the Gelfand-Fuchs 2-cocycle above and is known under the same
name. The modified Gelfand-Fuchs cocycle is $\sltwo$--invariant, where
$\sltwo$ is a subalgebra of $\vir^{\Bbb C}$ generated by $e_{-1}$, $e_0$ and
$e_1$, so it is handier in practice. Below we shall use the mo\-di\-fied
version of Gelfand-Fuchs 2-cocycle under this name only. In the irreducible
representation the central element $c$ of the Virasoro algebra is mapped to
a scalar operator, which is proportional to the identity operator with a
coefficient called {\it the central charge}.

\head 2. $\HS$-- and $\TC$--projective representations of the Witt
algebra\endhead

\subhead 2.1. $\frak A$--projective representations [6]\endsubhead

\definition{Definition 1 [6]}

{\bf A} Let $\frak A$ be an associative algebra represented by linear
operators in the linear space $H$ and $\frak g$ be the Lie algebra. The
linear mapping $T:\frak g\mapsto\End(H)$ is called {\it a
$\frak A$--projective representation} iff for any $X$ and $Y$ from $\frak g$
there exists an element of the algebra $\frak A$ represented by the operator
$A_{XY}$ such that $$[T(X),T(Y)]-T([X,Y])=A_{XY}.$$
The operators $A_{XY}$ are called {\it deviations\/} for the
$\frak A$--projective representation $T$.

{\bf B} Let $\frak A$ be an associative algebra with an
involution $*$ symmetrically represented in the Hilbert space $H$. If
$\frak g$ is a Lie algebra with an involution $*$ then its
$\frak A$--projective representation $T$ in the space $H$ is called
{\it symmetric\/} iff for all elements $a$ from $\frak g$ $T(a^*)=T^*(a)$.
Let $\frak g$ be $\Bbb Z$--graded Lie algebra
($\frak g=\oplus_{n\in\Bbb Z}\frak g_n$) with an involution $*$ such that
$\frak g^*_n=\frak g_{-n}$ and involution is identical on the subalgebra
$\frak g_0$. Let us extend the $\Bbb Z$--grading and the involution $*$
from the Lie algebra $\frak g$ to the tensor algebra $\bold
T^{\cdot}(\frak g)$. A symmetric $\frak A$--projective representation
of $\frak g$ ia called {\it absolutely symmetric\/} iff for any element
$a$ of the algebra $\bold T^{\cdot}(\frak g)$ such that $\deg(a)=0$
the identity $T(a)=T^*(a)$ holds (here the representation $T$ of the
Lie algebra $\frak g$ in $H$ is extended to the mapping of $\bold
T^{\cdot}(\frak g)$ into $\End(H)$.

{\bf C} A $\frak A$--projective representation $T$ of the
Lie algebra $\frak g$ in the linear space $H$ is called {\it almost
absolutely closed\/} iff for any natural $n$ and arbitrary elements
$X_0, X_1, X_2,\ldots X_{n+1}$ of the Lie algebra $\frak g$ there exists
an element $\varphi(X_0,X_1,X_2,\ldots\mathbreak X_{n+1})$ of the algebra
$\frak g$ such that
$$[\,\ldots[[T(X_0),T(X_1)],T(X_2)],\ldots\!,T(X_{n+1})]\!\equiv\!
T(\varphi(X_0,X_1,X_2,\ldots X_{n+1}))\!\!\!\!\pmod{\frak A},$$
here $\frak A$ is considered as represented in $\End(H)$. An almost
absolutely closed $\frak A$--projective representation $T$ of the Lie algebra
$\frak g$ in the linear space $H$ is called {\it absolutely closed\/} iff
$\varphi(\cdot,\ldots,\cdot)\equiv 0$.
\enddefinition

If $H$ is infinite dimensional then the representation may be realized
by the unbounded operators with suitable domains of definition.

\remark{Remark 1} The standard projective representation is a particular
case of Def.1A if the algebra $\frak A$ is one-dimensional and acts in
$H$ by scalar operators.
\endremark

\remark{Remark 2} If $H$ is a Hilbert (or pre-Hilbert) space then
one may consider the algebra $\HS$ of all Hilbert-Schmidt operators
as $\frak A$. It is possible to consider the algebras $\bnd$, $\K$, $\TC$
and $\FR$ of all bounded, compact, trace class and finite-rank operators,
too.
\endremark

\remark{Remark 3} One may define $\frak a$--projective representations
of Lie algebra $\frak g$, where $\frak a$ is a Lie algebra in the similar
way. Any $\frak A$--projective representation, where $\frak A$ is an
associative algebra, is a $\frak A_{[\cdot,\cdot]}$--projective
representation, where $\frak A_{[\cdot,\cdot]}$ is a commutator Lie algebra
of the associative algebra $\frak A$.
\endremark

\subhead 2.2. $\HS$-- and $\TC$--projective representations of Lie algebras
\endsubhead
Below we shall consider either $\HS$-- or $\TC$--projective absolutely
closed representations of Lie algebras only. Let's introduce some useful
concepts.

\definition{Definition 2A} Let $T$ be a $\TC$--projective representation
of Lie algebra $\frak g$. Then the {\it fundamental form\/} $\alpha_{XY}$
on $\frak g$ for the representation $T$ will be defined as
$$\alpha_{XY}=\tr A_{XY},$$
where $A_{XY}$ are deviations for representation $T$, $X,Y\in\frak g$.
\enddefinition

Let us denote by $\Ker A=\{X\in\frak g:(\forall Y\in\frak g) A_{XY}=0\}$
and $\Ker\alpha=\{X\in\frak g:(\forall Y\in\frak g) \alpha_{XY}=0\}$;
$\frak g_0=\Ker A$ is a Lie subalgebra of $\frak g$. Note that the
fundamental form $\alpha$ is $\frak g_0$--invariant. Let us suppose that
$\Ker A=\Ker\alpha$ then one can define the formal series
$$A=\alpha^{XY}A_{XY}$$
(here and below, the Einstein rule of summation is assumed). If it converges
absolutely in certain topology in some class $\Cal O$ of operators, the
result will be called {\it the mean deviation}.

\remark{Remark 4} The operator $A$ is $\frak g_0$--invariant.
\endremark

If $T$ is (operator) irreducible representation of $\frak g_0$ then one can
define the {\it de\-viation constant\/} $\alpha$ as
$$A=\alpha E,$$
where $E$ is the identical operator (note that $\alpha$ is not obligatory
equal to $1$ in the considered infinite--dimensional case).

Let us also define $B_{XYZ}=[T(X),A_{YZ}]$ and $\beta_{XYZ}=\tr B_{XYZ}$
(the least expressions exist because $T$ is absolutely closed).

\remark{Remark 5}
$$d\alpha_{XYZ}+(\beta_{XYZ}+\beta_{YZX}+\beta_{ZXY})=0,$$
where
$$d\alpha_{XYZ}=\alpha_{X,[Y,Z]}+\alpha_{Y,[Z,X]}+\alpha_{Z,[X,Y]}.$$
\endremark

The equalities follow from the identity
$$(A_{X,[Y,Z]}+B_{XYZ})+(A_{Y,[Z,X]}+B_{YZX})+(A_{Z,[X,Y]}+B_{ZXY})=0.$$

\definition{Definition 2B} Let $T$ be a $\HS$--projective representation
of Lie algebra $\frak g$. Then the {\it curvature form\/} $R\in S^2\Lambda^2
(\frak g^*)$ of $T$ is defined as
$$R_{XY,UV}=\tr(A_{XY}A_{UV}).$$
If $T$ is also $\TC$--projective representation and $\Ker A=\ker\alpha$ then
let us introduce the {\it curvature operator\/}
$$R_{XY,U}^V=\alpha^{ZV}R_{XY,UZ}.$$
If the series $R^V_{XV,Y}$ and $R^V_{XY,V}$ absolutely converges they
will be called the {\it Ricci tensors\/} of $\TC$--projective representation
$T$ of the first and the second kind, res\-pec\-ti\-ve\-ly, and denoted by
$R^{(i)}_{XY}$ ($i=1,2$). The $\TC$--projective representation $T$ will be
called {\it einsteinian\/} (of the 1st or 2nd kind) iff $R^{(i)}_{XY}$ is
propostional to the fundamental form $\alpha$:
$$R^{(i)}_{XY}=\varrho^{(i)}\alpha_{XY}.$$
\enddefinition

Such definition of Ricci tensors is motivated by analogies from K\"ahler
geometry, where Ricci tensor may be defined equivalently via curvature
operator as $R^V_{XV,Y}$ or as $R^V_{XY,V}$ [15].

\remark{Remark 6} Let $T$ be a $\TC$--projective representation of
Lie algebra $\frak g$, which is an (operator) irreducible representation of
$\frak g_0$ (and $\Ker A=\Ker\alpha$), and the formal series for the mean
deviation $A$ absolutely converges then the series for $R^{(2)}_{XY}$ also
converges and $\varrho^{(2)}=\alpha$.
\endremark

\subhead 2.3. $\HS$--projective representation of the Witt algebra by
$q_R$--conformal symmetries in the unitarizable Verma modules $V_h$ over
the Lie algebra $\sltwo$
\endsubhead
Note that the Witt algebra $\frak w^{\Bbb C}$ admits a natural involution $*$.

\proclaim{Proposition 3 [6]}\it The generators $L_k$ ($k\in\Bbb Z$) of
$q_R$--conformal symmetries in the unitarizable Verma module $V_h$ over the
Lie algebra $\sltwo$ realize an absolutely symmetric $\HS$--projective
representation of the Witt algebra $\frak w^{\Bbb C}$.
\endproclaim

\remark{Remark 7} The $\HS$--projective representations of the Witt algebra
in the uni\-ta\-ri\-zable Verma modules over $\sltwo$ are absolutely closed.
\endremark

\remark{Remark 8} Proposition 3 is generalized for the pseudo-unitary case
with the sub\-s\-ti\-tution of the class $\HS$ of Hilbert-Schmidt operators
to the class $\K$ of compact operators.
\endremark

\proclaim{Theorem 1} The $\HS$--projective representations of the Witt
algebra $\frak w^{\Bbb C}$ in the unitarizable Verma modules $V_h$ over the
Lie algebra $\sltwo$ by $q_R$--conformal sym\-met\-ries are $\TC$--projective
representations. The fundamental form $\alpha_{XY}$ in the basis $e_k$
has the form
$$\alpha_{ij}=\frac{6h^2-6h+1}6(j^3-j)\delta(i+j).$$
\endproclaim

\demo{Proof} In the considered case $\Ker A=\sltwo$, hence the fundamental
form $\alpha_{XY}$ is $\sltwo$--invariant. All $\sltwo$--invariant
2-forms on the Witt algebra $\frak w^{\Bbb C}$ are proportional to the
Gelfand--Fuchs cocycle. To determine the multiplier it is suf\-ficient to
calculate $\alpha_{2,-2}$ using explicit formulas for $q_R$--conformal
symmetries.
\enddemo

\remark{Remark 9} If $h=\frac12\pm\frac{\sqrt{3}}6$ then $\alpha_{XY}=0$.
\endremark

\remark{Remark 10} The fundamental form $\alpha_{XY}$ is closed, i.e.
$d\alpha_{XYZ}=0$, moreover $\beta_{XYZ}=0$.
\endremark

The corollary below follows from Theorem 1 and Remark 10.

\proclaim{Corollary} The deviations $A_{XY}$ define an 1-dimensional
central extension of $\frak w^{\Bbb C}$ by $\im A/(\im A\cap\TC_0)$
(here $\im A\subseteq\TC$ is the image of the mapping
$A:\Lambda^2(\frak w^{\Bbb C})\mapsto\TC$ and $\TC_0$ is
the class of all operators from $\TC$ of zero trace) with the central
charge
$$c=-2(6h^2-6h+1).$$
Otherwords, the $q_R$--conformal symmetries define the $\TC_0$--projective
representation ($\TC_0$ is a Lie algebra) of the Virasoro algebra
$\vir^{\Bbb C}$ with this central charge.
\endproclaim

\proclaim{Theorem 2A} For the $\TC$--projective representations of
the Witt algebra in the unitarizable Verma modules $V_h$ ($h\ne\frac12\pm
\frac{\sqrt{3}}6$) over $\sltwo$ the kernels $\Ker A$ and $\Ker\alpha$
coincide, therefore, the formal series $\alpha^{XY}A_{XY}$ are defined.
The series absolutely converges if $h>\frac12$, $h\ne\frac12+\frac{\sqrt{3}}6$.
The unitarizable Verma modules $V_h$ are $\sltwo$--irreducible, therefore,
the mean deviation $A$ is proportional to the identical operator $E$.
Precisely,
$$\alpha=-\frac3{2h-1}.$$
\endproclaim

An estimation of the convergence and formulas for $\alpha$ follow from
the explicit computations.

\proclaim{Corollary} The $\TC$--projective representations of the Witt
algebra in the uni\-ta\-ri\-zable Verma modules $V_h$ ($h>\frac12$,
$h\ne\frac12+\frac{\sqrt{3}}6$) over $\sltwo$ are einsteinian of the second
kind with $\varrho^{(2)}=\alpha$.
\endproclaim

\remark{Remark 11} The central charge $c$ and the deviation constant $\alpha$
obey the identity
$$\alpha^2(1-c)=27.$$
\endremark

\proclaim{Theorem 2B} For $h>\frac12$, $h\ne\frac12+\frac{\sqrt{3}}6$
the formal series for the Ricci tensor $R^{(1)}_{XY}$ converges and the
$\TC$--projective representations of the Witt algebra $\frak w^{\Bbb C}$
in the Verma modules $V_h$ over $\sltwo$ are einsteinian of the first kind.
\endproclaim

\demo{Proof} An estimation of convergence follows from the explicit
formulas for $q_R$--conformal symmetries. If the Ricci tensor $R^{(1)}_{XY}$
exists then it is $\sltwo$--invariant. However, any $\sltwo$--invariant
2-form on the Witt algebra $\frak w^{\Bbb C}$ is proportional to the
Gelfand-Fuchs cocycle.
\enddemo

\remark{Remark 12A} It is not known whether the Ricci constants $\varrho^{(1)}$
and $\varrho^{(2)}$ coincide.
\endremark

\remark{Remark 12B} It is not known also whether the Ricci identity [15,16]
$$R_{XY,ZV}+R_{YZ,XV}+R_{ZX,YV}=0$$
holds for the curvature form $R$.
\endremark

It is rather interesting to calculate some ``differential geometric''
characteristics of $q_R$--conformal symmetries such as ``scalar curvature''
$$K_{ij}=\frac{R_{i,-j,j,-i}}{\alpha_{i,-i}\alpha_{j,-j}}$$
and to verify whether the other identities of Riemann geometry besides Ricci
one [15,16] hold for the curvature form $R$.

\remark{Remark 13} The constants $c$ and $\alpha$ may be written in terms of
$q_R$ instead of $h$ as
$$c=1-3q_R^{-2},\qquad\qquad\alpha=-3q_R.$$
Hence,
$$\alpha\underset{q_R\to0}\to\longrightarrow0.$$
The central charge $c$ is equal to $0$ iff $q_R=\dfrac1{\sqrt{3}}$.
\endremark

\proclaim{Theorem 2C} Let put
$$B_X=\alpha^{YZ}B_{XYZ}=\alpha^{YZ}A_{[XY],Z}$$
for $h>\frac12$, $h\ne\frac12+\frac{\sqrt{3}}6$. Then
$$B_X=\alpha T(X).$$
\endproclaim

\demo{Proof} Note that the operators $B_X$ form the family of tensor
operators in the Verma module $V_h$, which transforms under $\sltwo$ just
as $e_k$, and, hence, $B_X$ is proportional to $T(X)$ because such family
of tensor operators is unique in $V_h$. To calculate the multiplier one
should put $X=e_0$ and to apply both $B_X$ and $T(X)$ to the
extremal vector.
\enddemo

\remark{Remark 14} The construction of Theorem 3C allows to define a mapping
from $\End(\frak w^{\Bbb C})$ to $\End(V_h)$, namely, the operator
$Q\in\End(\frak w^{\Bbb C})$ is mapped to the operator
$$\hat Q_h=\frac1{\alpha}\alpha^{YZ}A_{Q(Y),Z}$$
in $V_h$ (certainly, the convergence of series depends on $Q$).
\endremark

It is very interesting to investigate the asymptotical properties
of the mapping $Q\mapsto\hat Q_h$ (as a mapping of one associative
algebra into another) at $h\to\frac12$, $h\to 1$ (cf.[6]) and $h\to\infty$.

\head 3. Infinite dimensional geometry of homogeneous K\"ahler manifolds
$M=\Diff_+(\Bbb S^1)/\Bbb S^1$ and $M_1=\Diff_+(\Bbb S^1)/\PSL(2,\Bbb R)$
\endhead

\subhead 3.1. Group $\Diff_+(\Bbb S^1)$ of diffeomorphisms of a circle
$\Bbb S^1$ and the Virasoro-Bott group $\Vir$
\endsubhead
Let $\Diff(\Bbb S^1)$ denote the group of all diffeomorphisms of the circle
$\Bbb S^1)$. The group manifold $\Diff(\Bbb S^1)$ splits into two
connected components, the subgroup $\Diff_+(\Bbb S^1)$ and the coset
$\Diff_-(\Bbb S^1)$. The diffeomorphisms in $\Diff_+(\Bbb S^1)$ preserve
the orientation on the circle $\Bbb S^1)$ and those in $\Diff_-(\Bbb S^1)$
reverse it. The Lie algebra of $\Diff_+(\Bbb S^1)$ can be identified with
$\Vect(\Bbb S^1)$.

The infinite-dimensional group $\Vir$ corresponding to the Virasoro algebra
$\vir$ (more precisely, to the central extension of the Lie algebra
$\Vect(\Bbb S^1)$ defined by the Gelfand-Fuchs cocycle, whereas the
Virasoro algebra $\vir$ is an extension of the real form $\frak w$ of the
Witt algebra $\frak w^{\Bbb C}$) is a central extension of the group
$\Diff(\Bbb S^1)$. The corresponding 2-cocycle was calculated by R.Bott [17]:
$$c(g_1,g_2)=\int\log(g_1'\circ g_2)\,\log(g_2').$$
The group $\Vir$ is called the Virasoro-Bott group.

\subhead 3.2. Flag manifold $M=\Diff_+(\Bbb S^1)/\Bbb S^1$ of the
Virasoro-Bott group\endsubhead
The flag manifold $M$ of the Virasoro-Bott group is a homogeneous space
with transformation group $\Diff_+(\Bbb S^1)$ and isotropy group $\Bbb S^1)$.
There exist several different realizations of this manifold [18-20].

{\it Algebraic realization}. The space $M$ can be realized as a conjugacy
class in the group $\Diff_+(\Bbb S^1)$ or in the Virasoro-Bott group $\Vir$.

{\it Probabilistic realization}. Let $P$ be the space of real probability
measures $\mu=u(t)\,dt$ with smooth positive density $u(t)$ on $\Bbb S^1$.
The group $\Diff_+(\Bbb S^1)$ naturally acts on $P$ by the formula
$$g:u(t)\,dt\mapsto u(g^{-1}(t))\,dg^{-1}(t).$$
The action is transitive and the stabilizer of the point $(2\pi)^{-1}\,dt$
is isomorphic to $\Bbb S^1$, therefore, $P$ can be identified with $M$.

{\it Orbital realization}. The space $M$ can be considered as an orbit of the
coadjoint representation of $\Diff_+(\Bbb S^1)$ or $\Vir$. Namely, the
elements of the dual space $\vir^*$ of the Virasoro algebra $\vir$ are
identified with the pairs $(p(t)\,dt^2,b)$; the coadjoint action of $\Vir$
has the form
$$K(g)(p,b)=(gp-bS(g),b),$$
where
$$S(g)=\frac{g'''}{g'}-\frac32\left(\frac{g''}{g'}\right)^2$$
is the Schwarzian (the Schwarz derivative) and $gp$ denotes the natural action
of $g$ on the quadratic differential $p$. The orbit of the point $(a\cdot
dt^2,b)$ coincides with $M$ if $a/b\ne-n^2/2$, $n=1,2,3,\cdots$. Therefore,
a family $\omega_{a,b}$ of symplectic structures (Kirillov forms) is defined
on $M$ (cf.[21]).

{\it Analytic realization}. Let us consider the space $S$ of univalent
functions on the unit disk $D_+$ [22-24]. The Taylor coefficients $c_1, c_2,
c_3,\cdots$ in the expansion
$$f(z)=z+c_1z^2+c_2z^3+\cdots+c_nz^{n+1}+\cdots$$
form a coordinate system on $S$. The class $S$ can be naturally identified
with $M$ via the Kirillov construction [20]. The Lie algebra $\Vect(\Bbb S^1)$
acts on $S$ by the formulas
$$L_v(f(z))=if^2(z)\oint\left(\frac{wf'(w)}{f(w)}\right)^2\frac{v(w)}{f(w)-
f(z)}\frac{dw}w.$$
The Kirillov construction supply $M$ by the complex structure. The symplectic
structure $\omega_{a,b}$ coupled with the complex structure determines a
K\"ahler metric $w_{a,b}$ on $M$. More detailed information on the
infinite-dimensional geometry of the flag manifold $M$ is contained in [2]
(see also refs wherein). Note only that the curvature tensors of the K\"ahler
connections on $M$ were calculated in [1].

\subhead 3.3. Infinite dimensional K\"ahler manifold $M_1=\Diff_+(\Bbb
S^1)/\PSL(2,\Bbb R)$\endsubhead
The subgroup $\Bbb S^1$ is contained in each of the subgroups $H_k$,
$k=1,2,3,\cdots$, generated by the generators $ie_0$, $e_k-e_{-k}$ and
$i(e_k+e_{-k})$. The subgroup $H_k$ is isomorphic to the $k$-folded
covering of the group $H_1=\PSL(2,\Bbb R)$ and acts on $\Bbb S^1$ by the
formulas
$$z\mapsto\left(\frac{az^k+b}{\bar bz^k+\bar a}\right)^{1/k}=
\alpha z\left(\frac{1+bz^{-k}}{1+\bar bz^k}\right)^{1/k},$$
where $\alpha=(a/|a|)^{2/k}$ is an univalued function on $H_k$.

The homogeneous space $M_k=\Diff_+(\Bbb S^1)/H_k$ is a symplectic manifold
and can be identified with the orbit with $a/b=-k^2/2$ in the coadjoint
representation of the Virasoro-Bott group. The manifold $M_k$ is a quotient
of $M$ by the $\Diff_+(\Bbb S^1)$-invariant foliation $\Cal F_k$ generated
by $s_k$ and $c_k$ (considered as elements of a reductive basis on $M$ [15]).
Note that though $M$ is a reductive space the manifolds $M_k$ are not
reductive. The foliations $\Cal F_k$ are complex foliations, which were
discovered by V.Yu.Ovsienko and O.D.Ovsienko. The foliation $\Cal F_1$
is holomorphic whereas other foliations are not holomorphic. Hence,
the almost complex structure on $M_1$ is integrable and, therefore,
$M_1=\Diff_+(\Bbb S^1)/\PSL(2,\Bbb R)$ is a homogeneous K\"ahler manifold.
This is not true for $k\ge2$.

The main characteristics of the K\"ahler metric on $M_1$ were calculated in
[1] (really, there were calculated their back-liftings to $M$). In
particular, these K\"ahler metrics are einsteinian.

\subhead 3.4. Deformed holomorphic tangent bundles on
$M_1=\Diff_+(\Bbb S^1)/\PSL(2,\Bbb R)$, $\Diff_+(\Bbb S^1)$-in\-va\-riant
hermitean connections and their curvatures\endsubhead
Note that one can deform the holomorphic tangent bundle $T^{1,0}(M_1)$ on
$M_1$ in the following way. Let us consider the principal
$\PSL(2,\Bbb R)$--bundle $\Cal P$ over the homogeneous manifold $M_1$
and the associated holomorphic bundle $T_h(M_1)$ with the fiber iso\-mor\-phic
to $V_h$. The bundle $T_h(M_1)$ is supplied by the holomorphic action of
the Witt algebra $\frak w^{\Bbb C}$, which is integrable to the action of
the group $\Diff_+(\Bbb S^1)$. The holomorphic tangent bundle $T^{1,0}(M_1)$
is just the bundle $T_2(M_1)$. The hermitean form in $V_h$ ($h>\frac12$)
supplies $T_h(M_1)$ by a structure of the hermitean bundle.

Now we are able to formulate the main theorem of the article.

\proclaim{Theorem 3A} There exists a unique holomorphic $\frak w^{\Bbb
C}$--invariant ($\Diff_+(\Bbb S^1)$--in\-va\-riant) hermitean connection in
$T_h(M_1)$. The Nomizu operators of the connection coincides with
the generators of $q_R$--conformal symmetries $L_k$, whereas its curvature
coincides with the deviations $A_{XY}$.
\endproclaim

Note that though $M_1$ is not reductive the Nomizu operators are correctly
defined (one can use the back-lifting to the reductive space $M$).

\demo{Proof} One should use the definition of $q_R$--conformal symmetries.
Indeed, $L_k$ ($k\le1$) define the natural action of the corresponding
generators of the Witt algebra in the fiber $V_h$ of the bundle $T_h(M_1)$,
whereas $L_{-k}=L_k^*$.
\enddemo

\proclaim{Corollary} The deviations $A_{XY}$ coincides with the curvature
of K\"ahler metric on $M_1$ [1] at $h=2$.
\endproclaim

\remark{Remark 15} The curvature form $R$ of the $\HS$--projective
representation of the Witt algebra $\frak w^{\Bbb C}$ in the Verma module
$V_h$ over the Lie algebra $\sltwo$ is just the trace of a square of
the curvature of holomorphic $\frak w^{\Bbb C}$--invariant hermitean
connection in the deformed holomorphic tangent bundle $T_h(M_1)$ over
the infinite dimensional K\"ahler manifold $M_1=\Diff_+(\Bbb
S^1)/\PSL(2,\Bbb R)$.
\endremark

Let us reformulate Theorem 1 in the following way.

\proclaim{Theorem 3B} The holomorphic $\frak w^{\Bbb C}$--invariant
hermitean connection in the de\-formed holomorphic tangent bundle
$T_h(M_1)$ over the infinite-dimensional K\"ahler manifold
$M_1=\Diff_+(\Bbb S^1)/\PSL(2,\Bbb R)$ defines the prequantization
connection [25] in the deteminant bundle $\det T_h(M_1)$ over this manifold,
which supplies $\det T_h(M_1)$ by an action of the Virasoro algebra
$\vir^{\Bbb C}$ with the central charge $$c=-2(6h^2-6h+1).$$
\endproclaim

\remark{Remark 16} At $h=2$ the central charge is equal to $-26$ as it was
once calculated in [26].
\endremark

\head Conclusions\endhead

Thus, we established a geometric interpretation of approximate
($\HS$--projective or $\TC$--projective) representations of the Witt algebra
$\frak w^{\Bbb C}$ by $q_R$--conformal sym\-met\-ries in the Verma modules
$V_h$ over the Lie algebra $\sltwo$ and
calculated some their characteristics. The generators of representations
coincide with the Nomizu ope\-ra\-tors of ho\-lo\-mor\-phic $\frak w^{\Bbb
C}$--invariant hermitean connections in the deformed ho\-lo\-mor\-phic tangent
bundles $T_h(M_1)$ over the infinite-dimensional K\"ahler manifold
$M_1=\Diff_+(\Bbb S^1)/\PSL(2,\Bbb R)$, whereas the deviations $A_{XY}$
of the approximate rep\-re\-sen\-tations coincide with the curvature operators
for these connections, which supply the determinant bundle $\det T_h(M_1)$
by a structure of the prequantization bundle over $M_1$. At $h=2$ the
geometric picture reduces to one considered in [1] (without any relation
to approximate representations) for ordinary tangent bundles.

\Refs
\roster
\item"[1]" Kirillov A.A., Juriev D.V., The K\"ahler geometry of the infinite
dimensional homogeneous space $M=\Diff_+(\Bbb S^1)/\Rot(\Bbb S^1)$.
Funct.~Anal.~Appl. 21(4) (1987) 284-293.
\item"[2]" Juriev D.V., The vocabulary of geometry and harmonic analysis on
the infinite-di\-men\-sional manifold $\Diff_+(\Bbb S^1)/\Bbb S^1$.
Adv.~Soviet Math. 2 (1991) 233-247.
\item"[3]" Juriev D., Quantum conformal field theory as infinite dimensional
noncommutative geo\-met\-ry. Russian Math.~Surveys 46(4) (1991) 135-163.
\item"[4]" Juriev D., Infinite dimensional geometry and quantum field theory
of strings. I-III. Alg.~Groups Geom. 11(2) (1994) 145-179, Russian
J.~Math.~Phys. 4(3) (1996) 287-314, J.~Geom.~Phys. 16 (1995) 275-300.
\item"[5]" Juriev D., String field theory and quantum groups. I.:
q-alg/9708009.
\item"[6]" Juriev D., Approximate representations and Virasoro algebra:
math.RT/9805001.
\item"[7]" Zhelobenko D.P., Representations of the reductive Lie algebras
[in Russian]. Moscow, Nauka, 1993.
\item"[8]" Dixmier J., Alg\`ebres enveloppantes. Gauthier-Villars, Paris,
1974.
\item"[9]" Barut A., Raczka R., Theory of group representations and
applications. PWN -- Polish Sci.~Publ. Warszawa, 1977.
\item"[10]" Juriev D.V., Complex projective geometry and quantum projective
field theory. Theor. Math.~Phys. 101(3) (1994) 1387-1403.
\item"[11]" Berezin F.A., Quantization in complex symmetric spaces [in
Russian]. Izvestiya AN SSSR. Ser.~matem. 39(2) (1975) 363-402.
\item"[12]" Klimek S., Lesniewski A., Quantum Riemann surfaces. I.
Commun.~Math.~Phys. 146 (1992) 103-122.
\item"[13]" Biedenharn L., Louck J., Angular momentum in quantum mechanics.
Theory and appli\-cations. Encycl.~Math.~Appl. V.8. Addison Wesley Publ.~Comp.
1981.
\item"[14]" Fuchs D.B., Cohomology of infinite-dimensional Lie algebras [in
Russian]. Moscow, Nau\-ka, 1984.
\item"[15]" Kobayashi Sh, Nomizu K., Foundations of differential geometry.
Intersci.~Publ., 1963/69.
\item"[16]" Rashevski{\v\i} P.K., Riemann geometry and tensor analysis
[in Russian]. Moscow, 1967.
\item"[17]" Bott R., The characteristic classes of groups of diffeomorphisms.
Enseign.~Math. 23 (1977) 209-220.
\item"[18]" Kirillov A.A., Infinite-dimensional Lie groups, their invariants
and representations. Lect. Notes Math. 970 (1982) 101-123.
\item"[19]" Segal G., Unitary representations of some infinite-dimensional
groups. Commun.~Math. Phys. 80 (1981) 301-342.
\item"[20]" Kirillov A.A., A K\"ahler structure on K-orbits of the group of
diffeomorphisms of the circle. Funct.Anal.Appl. 21 (1987) 322-324.
\item"[21]" Kirillov A.A., Elements of representation theory. Springer, 1976.
\item"[22]" Goluzin G.M., Geometric theory of functions of complex variables.
Amer.~Math.~Soc. 1968.
\item"[23]" Duren P.L., Univalent functions. Springer, 1983.
\item"[24]" Lehto O., Univalent functions and Teichm\"uller spaces. Springer,
1986.
\item"[25]" Kostant B., Quantization and unitary representations.
Lect.~Notes Math. 170 (1970) 87-208.
\item"[26]" Bowick M.J., Rajeev S.G., The holomorphic geometry of closed
bosonic strings and $\Diff(S^1)/S^1$. Nucl. Phys. B293 (1987) 348-384.
\endroster
\endRefs
\enddocument